# Multi-Panel Kendall Plot in Light of an ROC Curve Analysis Applied to Measuring Dependence

ALBERT VEXLER*

*Department of Biostatistics, The State University of New York at Buffalo, USA*

*Corresponding author. Email: avexler@buffalo.edu

GEORGIOS AFENDRAS

*Department of Biostatistics and Jacobs School of Medicine and Biomedical Sciences,The State University of New York at Buffalo, USA*

MARIANTHI MARKATOU

*Department of Biostatistics and Jacobs School of Medicine and Biomedical Sciences,The State University of New York at Buffalo, USA*

The Kendall plot (K-plot) is a plot measuring dependence between the components of a bivariate random variable. The K-plot graphs the Kendall distribution function against the distribution function of $VU$, where $V$ and $U$ are independent uniform $[0,1]$ random variables. We associate K-plots with the receiver operating characteristic (ROC) curve, a well-accepted graphical tool in biostatistics for evaluating the ability of a biomarker to discriminate between two populations. The most commonly used global index of diagnostic accuracy of biomarkers is the area under the ROC curve (AUC). In parallel with the AUC, we propose a novel strategy to measure association between random variables from a continuous bivariate distribution. First, we discuss why the area under the conventional Kendall curve (AUK) cannot be used as an index of dependence. We then suggest a simple and meaningful extension of the definition of the K-plots, and define an index of dependence that is based on AUK. This measure characterizes a wide range of two-variable relationships, thereby completely detecting the underlying dependence structure. Properties of the proposed index satisfy the mathematical definition of a measure. Finally, simulations and real data examples illustrate the applicability of the proposed method.

## 1. Introduction

Assume we want to assess whether the random variables $X$ and $Y$ from a continuous bivariate distribution function $H$ with corresponding marginals $F$ and $G$ are stochastically independent or not. Various publications have proposed and studied different methods to measure a distance between the distributions $H(x,y)$ and $F(x)G(y)$ to quantify dependence for the pair $(X,Y)$ (e.g., Balakrishnan & Lai, 2009; Reshef et al., 2011; Vexler et al., 2014). Graphical display tools to evaluate a rich class of bivariate dependence structures (Fisher & Switzer, 1985, 2001; Bjerve & Doksum, 1993; Jones, 1996; Jones & Koch, 2003) include a K-plot procedure (Genest & Boies, 2003; Gargouri-Ellouze & Bargaoui, 2009). The K-plot suggests to plot the Kendall distribution function $K(t) = \Pr\{H(X,Y) < t\}$ against the function $\Pr\{F(X)G(Y) < t \mid X, Y \text{ are independent}\}$. This probability is equivalent to $\Pr(VU < t) = t - t\log(t)$ with $t \in [0,1]$, where the random variables $V$ and $U$ are independent and uniformly $[0,1]$ distributed (e.g., Genest & Rivest, 1993; Nelsen et al., 2003). The recent statistical literature has gradually recognized that the function $K(t)$ contains useful information regarding the dependence structure underlying $H(x,y)$. The visualization of the curve $(K(t), t - t\log(t))$ compared with the diagonal line $(t - t\log(t), t - t\log(t))$ in $[0,1] \times [0,1]$, for $t \in [0,1]$, provides a rich source of information about association related to pairs of a random sample from $H$. For example, interesting real data applications of the K-plots analysis are shown in Boero et al. (2011) and Eslamian (2014).

The graphical K-plot tool can be considered in light of the ROC curve methodology. Receiver operating characteristic curve analysis is a well-accepted statistical method for evaluating the discriminatory ability of biomarkers (e.g., Shapiro, 1999; Vexler et al., 2008). An ROC curve plots the true positive rates of a biomarker versus its false-positive rates for various thresholds of the test result. It is a convenient way of evaluating diagnostic biomarkers because the ROC curve places tests on the same scale where they can be compared for accuracy. Let the variables $Z$ and $W$ represent values of a biomarker associated with diseased ($Z$) and healthy ($W$) populations, respectively. The ROC curve is usually depicted by plotting the two-dimensional curve $(1 - F_Z(u), F_W(u))$, for $u \in (-\infty, \infty)$, that is a graph of $\text{ROC}(t) = 1 - F_Z\left\{F_W^{-1}(1-t)\right\}$, $t \in [0,1]$, where $F_Z$ and $F_W$ are the distribution functions of $Z$ and $W$, and $F^{-1}$ is the inverse function of $F$. The ROC curve displays distances between the distribution functions $F_Z$ and $F_W$, comparing the curve $(1 - F_Z(u), F_W(u))$ with the diagonal line $(1 - F_W(u), F_W(u))$ in $[0,1] \times [0,1]$, for $u \in (-\infty, \infty)$. The area under the ROC curve (AUC) is a common and well developed index that summarizes the information contained in the ROC curve. Bamber (1975) showed that

$$\text{AUC} = \int (1 - F_Z(u)) \mathrm{d}F_W(u) = \Pr(Z > W).$$

Obviously, the closer the AUC is to 1, the better the diagnostic accuracy of the biomarker in terms of the distances between the distribution functions $F_Z$ and $F_W$. The case $F_Z = F_W$ is reflected by AUC $= 1/2$ that corresponds to the area under the diagonal line. The AUC is also closely related to the Gini coefficient (Breiman et al., 1984), which is twice the area between the diagonal and the ROC curve.

By virtue of the generality of the graphical methods related to the Kendall and ROC procedures, it is interesting to investigate the possibility of developing a general index of dependence by employing the area under the Kendall curve (AUK),

$$\text{AUK} = \int_0^1 \Pr\{H(X,Y) < t\} \mathrm{d}\Pr(VU < t)$$
$$= \int_0^1 K(t) \mathrm{d}\{t - t\log(t)\} = -\int_0^1 K(t) \log(t) \mathrm{d}t.$$

In analogy to the AUC, the further the AUK is from $1/2$, the stronger the association between $X$ and $Y$. The case of AUK $= 1/2$, the area under the diagonal line, indicates a scenario in which independence $H(x,y) = F(x)G(y)$ is in effect.

In Sections 2, 3 we argue, based on axioms, that K-plots for measuring dependence need to be extended for more accurate dependence measurement. We prove that multi-panel K-plots should be considered to consistently evaluate dependence tasks. The simple extension of the K-plot approach implies an effective and meaningful AUK-based measure of dependence, a ramification of the broadly applicable AUC type tools. The proposed index aids to characterize a wide range of two-variable relationships via comprehensive detection of the underlying dependence structure. Note that, in the context of modern measures of dependence, the recently introduced Maximal Information Coefficient (MIC) provides similar scores to different dependence structures with equal noise (Reshef et al., 2011). In situations when the classical measures of dependence or MIC provide powerful outputs, the proposed approach is still meaningful as graphical concepts incorporating information regarding dependence structures into data analysis are not yet very well developed.

In Section 4, the theoretical properties of the new dependence measure are derived. We refer to the Appendix for technical derivations and proofs. The obtained results provide an easy way to compute and interpret the AUK as a measure of dependence. It turns out that the AUK-based measures preserve an ordering of "more associated" for bivariate distributions in the context shown in Schriever (1987). The AUK measure is countably sub-additive, satisfying the requirement to be a mathematical measure. The developed approach satisfies an affine invariance principle. To our knowledge, such properties are not proven for the recently introduced MIC measure.

We employ the Farlie-Gumbel-Morgenstern formula of the bivariate distribution function

$$H(x, y) = F(x)G(y)[1 + \gamma\{1 - F(x)\}\{1 - G(y)\}], \tag{1}$$

where $\gamma \in (-1, 1)$, to demonstrate an example with a closed-form solution for the AUK index that reflects a broad class of dependence structures (e.g., Schucany et al., 1978).

Furthermore, in Section 4 we outline the association between estimators of the AUKs and the likelihood principle. Section 5 presents a simulation study for a variety of distributions and sample sizes that illustrates the performance of our methods. Section 6 illustrates the applicability of the proposed approach through a real world example. A myocardial infarction disease data set, related to the accuracy of biomarkers to discriminate "disease" and "non-disease" populations, shows the proposed methodology can significantly outperform the MIC procedure in practice. In Section 7 we provide concluding remarks, rekindling interest in – and discussing the potential usefulness of – the proposed approach.

## 2. An issue related to the Kendall distribution function

In 1959 Rényi's fundamental article defined a set of axioms that a measure of dependence for a pair of random variables must satisfy. One of the axioms states that if the joint distribution of $(X, Y)$ is bivariate normal, with correlation coefficient $\rho$, then a measure of dependence should be a strictly increasing function of $|\rho|$. It is reasonable that, in this case, values of the measure at $\rho$ and $-\rho$ should be equal (Balakrishnan & Lai, 2009).

Consider, for example, the Kendall distribution function $K(t)$ when $H(x, y)$ corresponds to the density function

$$h(x, y) = \frac{1}{2\pi(1 - \rho^2)^{0.5}} \exp\left(-\frac{x^2 - 2\rho xy + y^2}{2(1 - \rho^2)}\right).$$

Figure 1 depicts the corresponding K-plots for different values of $\rho$. This figure clearly illustrates the concern that the plotted Kendall distribution function does not display information regarding the dependence structure of the bivariate normal distributions in an identical manner regarding $\rho$ and $-\rho$ with respect to the distances between the curves $(K(t), t - t\log(t))$ and the diagonal line $(t - t\log(t), t - t\log(t))$ in $[0, 1] \times [0, 1]$, for $t \in [0, 1]$. The area between the curve above the diagonal line and the diagonal line is always larger than the area between the diagonal line and the curve below it.

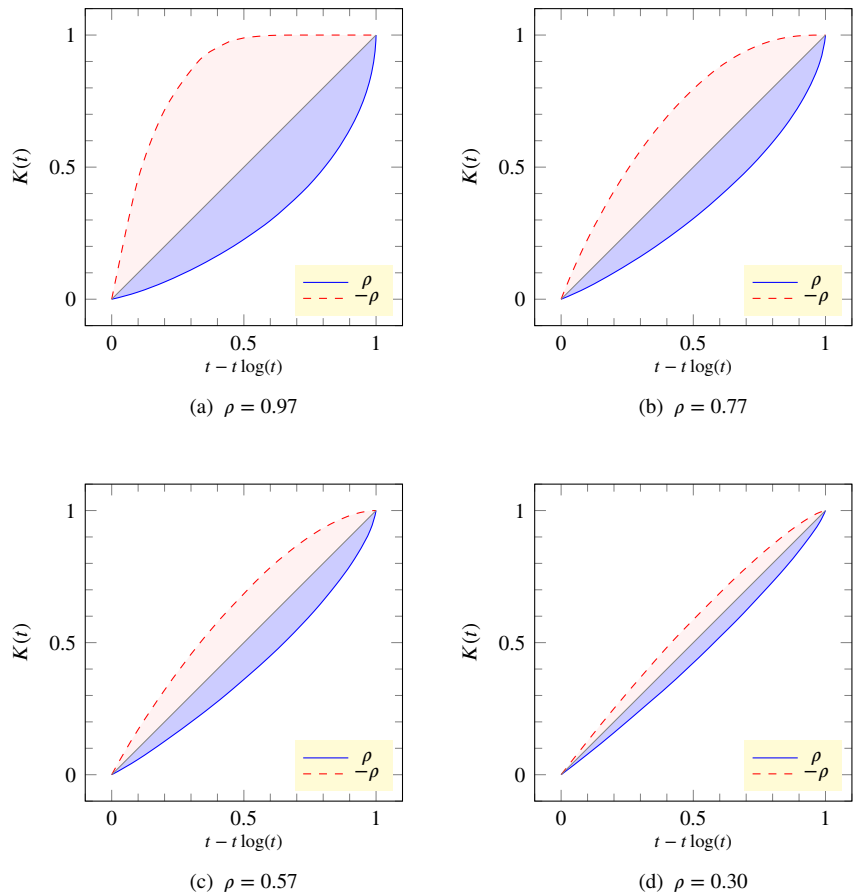

Figure 1: The $K(t)$ plotted against the function $t - t\log(t)$, for $t \in [0,1]$, when the joint distribution of $(X, Y)$ is a bivariate normal, with correlation $\rho$.

In this context, one can represent the bivariate normal random variable $(X, Y)$ in the form $(X, Y) = \left(X, \rho X + (1-\rho^2)^{1/2}\xi\right)$, where $X$ follows the standard normal distribution and $\xi$, independent of $X$, also follows the standard normal distribution. In this case, the distribution of $(X, Y)$ is $H(x, y) = \Pr\left\{X < x, \rho X + (1-\rho^2)^{1/2}\xi < y\right\}$ that obviously implies

$$\lim_{\rho \to 1} H(x, y) = \Pr(X < x, X < y) = \Phi(\min\{x, y\}),$$
$$\lim_{\rho \to -1} H(x, y) = \Pr(X < x, -X < y) = \max\{0, \Phi(x) - \Phi(-y)\},$$
$$\Phi(u) = (2\pi)^{-0.5} \int_{-\infty}^{u} \exp\left(-z^2/2\right) dz,$$

since $\lim_{\rho \to 1} Y = X$ and $\lim_{\rho \to -1} Y = -X$. Then the Kendall distribution function has different asymptotic shapes, $K(t) \to \Pr\{\Phi(X) < t\} = t$ when $\rho \to 1$ and $K(t) \to \Pr\{\Phi(X) - \Phi(X) < t\} = 1$ when $\rho \to -1$, for $t \in [0, 1]$. This conclusion also follows directly from the results shown in Kotz et al. (2000, Equations (46.38)–(46.45)). Thus, we obtain the following

$$\lim_{\rho \to 1} \text{AUK} = -\int_0^1 t \log(t) dt = 1/4, \quad \lim_{\rho \to -1} \text{AUK} = -\int_0^1 \log(t) dt = 1.$$

It is clear that the distance between AUK $= 1/2$ (the case of independence) and AUK $= 1/4$ (as $\rho \to 1$) and the distance between AUK $= 1/2$ and AUK $= 1$ (as $\rho \to -1$) are different, while those should be the same.

Furthermore, for illustrative purposes, we refer the reader to Figure 2 below that displays the corresponding AUK values depending on $\rho$. In Section 3, we also present an additional example that clearly illustrates the need of extending the K-plot.

Thus, the applied ROC/AUC technique assists to detect the aforementioned issue that motivates the need of extending the definition of the K-plot to a multi-panel K-plot in order to consistently evaluate dependence between the components of a bivariate random variable. The proposed extension of the K-plot based methodology will employ AUK type quantities. To this end, in Section 4, we derive essential properties of the AUK, in general.

## 3. The index of dependence based on multi-panel K-plots

In this section, we consider an extension of the K-plots that is more directly linked to the nature of dependence for pairs of random variables. Towards this end we define the probability function $H_1(x, y) = \Pr(X \geq x, Y < y)$. Assume the random variables $X$ and $Y$ are independent. Then, $H_1(x, y) = (1 - F(x))G(y)$. Since interest centers in measuring dependence, and $1 - U_1$, where $U_1 = F(X)$, is distributed according to $Unif[0,1]$ distribution, it is natural to compare the random variable $H_1(X,Y)$ with $VU$, $V$ and $U$ are independent $Unif[0,1]$ random variables. In a similar manner to the discussion presented above, one can define the probability functions $H_2(x,y) = \Pr(X < x, Y \geq y)$, $H_3(x,y) = \Pr(X \geq x, Y \geq y)$; and consider the Kendall distribution type functions $K_i(t) = \Pr\{H_i(X,Y) < t\}$, $i = 0, 1, 2, 3$, where $H_0(x,y) = H(x,y)$.

Accordingly, we extend the K-plot to multi-panel K-plot $(K_i(t), t - t\log(t))$, $i = 0, 1, 2, 3$. The multi-panel K-plot provides an informative approach to displaying positive and negative dependence structures of the bivariate distribution $H$ (see Section 4 for the corresponding results, in this context). Additionally, an associated measure of dependence is based on the vector $\boldsymbol{D}$ given as

$$\boldsymbol{D} = (\text{AUK}_0, \text{AUK}_1, \text{AUK}_2, \text{AUK}_3)^T, \quad \text{where } \text{AUK}_i = -\int_0^1 K_i(t)\log(t)\mathrm{d}t.$$

A strong argument in favor of the definition of the vector $\boldsymbol{D}$ is that the characterization of the association in a random sample from a continuous bivariate distribution is based on an affine invariance principle. Independence is preserved under affine transformations, hence it is natural to consider dependence measures that are affine invariant. Assume $a_X$, $a_Y$ are constant and $b_X$, $b_Y$ are arbitrary nonzero numbers. We define the transformed random variables $\widetilde{X} = a_X + b_X X$ and $\widetilde{Y} = a_Y + b_Y Y$ and consider their joint distribution function $\widetilde{H}$. It is clear that $\widetilde{H}(\widetilde{X}, \widetilde{Y}) = H(X,Y)$, when $b_X > 0$, $b_Y > 0$, whereas, e.g., $\widetilde{H}(\widetilde{X}, \widetilde{Y}) = \Pr\{X > (\widetilde{x}-a_X)/b_X, Y < (\widetilde{y}-a_Y)/b_Y\}$ and then $\widetilde{H}(\widetilde{X}, \widetilde{Y}) = H_1(X,Y)$, when $b_X < 0$, $b_Y > 0$. Thus, in this context, the four components of the vector $\boldsymbol{D}$, which are based on $H$, $H_1$, $H_2$ and $H_3$, should be considered to protect the affine invariance property of the measure.

An additional aspect to focus on the vector's $\boldsymbol{D}$ structure is a likelihood concept we discuss in Section 4.

Furthermore, the multi-panel K-plot yields a new global index of dependence based on the formula

$$I_{\text{AUK}} = (8/5)^{1/2}\|\boldsymbol{D} - \boldsymbol{\Delta}\| = (8/5)^{1/2}\left(\sum_{i=0}^{3}(\text{AUK}_i - 1/2)^2\right)^{1/2},$$

where $\boldsymbol{\Delta} = (1/2, 1/2, 1/2, 1/2)^T$ and the coefficient $(8/5)^{1/2}$ corresponds to the inverse of the Euclidean distance of the vector $\boldsymbol{D} - \boldsymbol{\Delta} = (-1/4, 1/2, 1/2, -1/4)^T$ computed under $(X,Y) = (X,X)$ (see Proposition 2 shown in Section 4). The case $H(x,y) = F(x)G(y)$ is reflected by $I_{\text{AUK}} = 0$.

Figure 2 shows the values of $I_{\text{AUK}}$ and AUK, when $H(x,y)$ is the bivariate normal distribution, plotted against the correlation coefficient $\rho \in (-1, 1)$. Notice that the AUK is not a symmetric function of $\rho$, while $I_{\text{AUK}}$ is.

In Section 4, we present detailed evaluations related to Rényi's axioms with respect to the proposed index. The structure of the vector $\boldsymbol{D}$ plays a key role in these considerations.

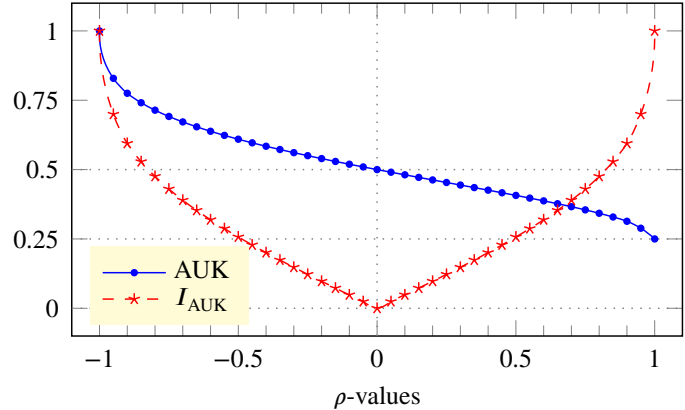

Figure 2: Values of the AUK and $I_{\text{AUK}}$ calculated for various values of $\rho \in (-1, 1)$.

In Section 4, we provide Proposition 3 to demonstrate that the AUK based measures are countably sub-additive.

Next we present a standardized index of dependence based on $I_{\text{AUK}}$ that ensures linear, one-to-one mapping with $|\rho|$ in the case of the bivariate normal distribution. Let $(X, Y)$ follow a bivariate normal distribution and the correlation between $X$ and $Y$ is $\rho$. Figure 2 shows that the $I_{\text{AUK}}$ is a monotonically increasing function, say $\eta$, of $|\rho|$ and ranges between $[0, 1]$; that is, $I_{\text{AUK}} = \eta(|\rho|)$, $|\rho| \in [0, 1]$. Therefore, there exists the inverse function of $\eta$, $\eta^{-1} \colon [0, 1] \to [0, 1]$. We define the *standardized* AUK-*index* by the relation

$$\text{st.}I_{\text{AUK}} = \eta^{-1}\left(I_{\text{AUK}}\right).$$

By definition, if $(X, Y)$ follows a bivariate normal distribution and the correlation between $X$ and $Y$ is $\rho$, then $\text{st.}I_{\text{AUK}} = |\rho|$. Using Lagrange interpolation method, we approximate numerically the function $\eta^{-1}$ by the polynomial

$$\tilde{\eta}^{-1}(t) = 2.070t + 0.061t^2 - 2.471t^3 + 1.307t^4 + 0.033t^5, \quad t \in [0, 1].$$

We denote the approximation of the standardized AUK-index by

$$\overline{I}_{\text{AUK}} = \tilde{\eta}^{-1}\left(I_{\text{AUK}}\right).$$

Figure 3 plots $I_{\text{AUK}}$ and $\overline{I}_{\text{AUK}}$, when $H(x,y)$ is the bivariate normal distribution, against $|\rho| \in [0,1]$. Observe that $\overline{I}_{\text{AUK}}$, as a function of $|\rho|$, is very close to the identity function.

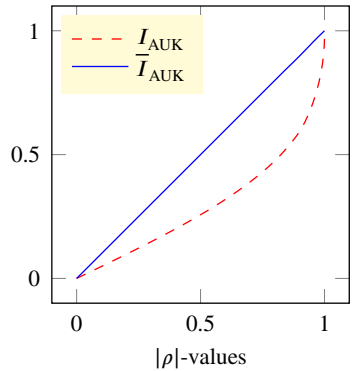

Figure 3: The plots of $I_{\text{AUK}}$ and $\overline{I}_{\text{AUK}}$ against $|\rho|$ in the bivariate normal case.

**Example:** Let $ABC$ be the triangle defined by the points $A(0,1)$, $B(-1,0)$ and $C(1,0)$. Suppose, for example, that the random vector $(X,Y)$ is uniformly distributed on the sides $AB$ and $AC$ of the triangle $ABC$. In this case, one can easily see that $K_0(t) = K_1(t) = 2t$ if $0 \leq t \leq 1/2$ and it equals 1 if $1/2 < t \leq 1$. Similarly, $K_2(t) = K_3(t) = 1/2 + t/2$ if $0 \leq t \leq 1/2$ and it equals 1 if $1/2 < t \leq 1$. (See the Appendix for details.) Thus, we have $\text{AUK}_0 = \text{AUK}_1 = 5/8 - \log(2)/4 \cong 0.4517$, $\text{AUK}_2 = \text{AUK}_3 = 25/32 - 3\log(2)/14 \cong 0.6513$ and $I_{\text{AUK}} = 0.284$, $\overline{I}_{\text{AUK}} = 0.545$. This also confirms the need for considering $\text{AUK}_i$, $i = 0,\ldots,3$ (not only one of the AUKs) to consistently display dependence between $X$ and $Y$, since, e.g., $|\text{AUK}_2 - 1/2| = |0.6513 - 0.5| > |0.4517 - 0.5| = |\text{AUK}_0 - 1/2|$. In this example, it is interesting to remark that $X \sim U[-1,1]$ and $Y = 1 - |X|$ implying $Cov(X,Y) = 0$. Note that, in a similar manner to the example above, one can show different scenarios in which $\text{AUK}_0 \neq AUK_1 \neq AUK_2 \neq AUK_3$.

## 4. Properties of AUK

In the ROC curve framework, Bamber (1975) expressed AUC in a simple form to facilitate computations based on the AUC. In a similar fashion, we present Proposition 1 that provides a simple formula for the AUK.

**Proposition 1.** *Suppose continuous random variables $X$ and $Y$ are distributed according to a bivariate distribution function $H$. Then,*

$$\text{AUK} = \Pr\{H(X,Y) < VU\} = \text{E}[1 - H(X,Y) + H(X,Y)\log\{H(X,Y)\}].$$

Proposition 1 is very useful both, for direct evaluation and for simulation studies based on the AUK.

Suppose that a random sample $(X_j, Y_j)$, $j = 1, \ldots, n$, has been drawn from a bivariate distribution $H$. By virtue of Proposition 1, we can estimate the AUK in a nonparametric manner via the statistic

$$\overline{\text{AUK}} = \frac{1}{n} \sum_{j=1}^{n} \left\{ 1 - \widehat{H}_j + \widehat{H}_j \log\left(\widehat{H}_j\right) \right\},$$

where $\widehat{H}_j = \widehat{H}(X_j, Y_j)$, replacing $H(\cdot, \cdot)$ by the corresponding empirical distribution $\widehat{H}(\cdot, \cdot)$. This estimator is much simpler than that we could obtain by using the direct definition of the AUK. It is clear that $\overline{\text{AUK}}$ is a consistent estimator of the AUK as $n \to \infty$. The statistic $\overline{\text{AUK}}$ has a structure similar to those of $U$-statistics. This can provide a wide range of evaluations of statistical properties related to $\overline{\text{AUK}}$ that deserves further strong empirical and methodological investigations.

In a similar manner to Proposition 1 we can express the quantities $\text{AUK}_i$, $i = 1, 2, 3$.

In Section 6 we employ nonparametric estimation of the AUK to perform straightforwardly resampling procedures for evaluating associations between biomarkers related to the myocardial infarction disease.

**The maximum likelihood point of view.** Consider the empirical estimation of the total $\text{AUK}_T = \sum_{i=0}^{3} \text{AUK}_i$ in the form

$$\frac{1}{n} \sum_{j=1}^{n} \left\{ 1 - \widehat{H}_j + \widehat{H}_j \log\left(\widehat{H}_j\right) \right\}$$

$$+ \frac{1}{n} \sum_{j=1}^{n} \left\{ 1 - \widehat{G}_j + \widehat{H}_j + \left(\widehat{G}_j - \widehat{H}_j\right) \log\left(\widehat{G}_j - \widehat{H}_j\right) \right\}$$

$$+ \frac{1}{n} \sum_{j=1}^{n} \left\{ 1 - \widehat{F}_j + \widehat{H}_j + \left(\widehat{F}_j - \widehat{H}_j\right) \log\left(\widehat{F}_j - \widehat{H}_j\right) \right\}$$

$$+ \frac{1}{n} \sum_{j=1}^{n} \left\{ \widehat{F}_j + \widehat{G}_j - \widehat{H}_j + \left(1 - \widehat{F}_j - \widehat{G}_j + \widehat{H}_j\right) \log\left(1 - \widehat{F}_j - \widehat{G}_j + \widehat{H}_j\right) \right\}$$

$$= 3 + \frac{1}{n} \sum_{j=1}^{n} \left\{ \widehat{H}_j \log\left(\widehat{H}_j\right) + \left(\widehat{G}_j - \widehat{H}_j\right) \log\left(\widehat{G}_j - \widehat{H}_j\right) + \left(\widehat{F}_j - \widehat{H}_j\right) \right.$$
$$\left. \log\left(\widehat{F}_j - \widehat{H}_j\right) + \left(1 - \widehat{F}_j - \widehat{G}_j + \widehat{H}_j\right) \log\left(1 - \widehat{F}_j - \widehat{G}_j + \widehat{H}_j\right) \right\},$$

where $\widehat{H}_j, \widehat{F}_j, \widehat{G}_j$ are the empirical estimators of $H(X_j, Y_j), F(X_j), G(Y_j)$ respectively. It is clear that this formal notation can be linked to the log maximum likelihood function based on the events $(X_j \leq X_{j'}, Y_j \leq Y_{j'})$, $(X_j \geq X_{j'}, Y_j \leq Y_{j'})$, $(X_j \leq X_{j'}, Y_j \geq Y_{j'})$ and $(X_j \geq X_{j'}, Y_j \geq Y_{j'})$.

**An ordering for dependence.** Schriever (1987) defined the "more associated" -ordering for bivariate distributions. By virtue of Proposition 1, we have $1 - \text{AUK} = \iint J\{H(x,y)\} dH(x,y)$, where the function $J(u) = u - u \log(u)$ increases and is upper convex for $u \in (0,1)$. Then Example 3.2 in Schriever (1987) can be directly adapted to show that the proposed AUK based measures preserve "more concordant" -ordering for dependence. We refer the reader to Schriever (1987) for details regarding the ordering for dependence. Consequently, the differences $(1 - \text{AUK}_0) - 1/2$, $\text{AUK}_1 - 1/2$, $\text{AUK}_2 - 1/2$ and $(1 - \text{AUK}_3) - 1/2$ measure positive dependence, whereas $\text{AUK}_0 - 1/2$, $(1 - \text{AUK}_1) - 1/2$, $(1 - \text{AUK}_2) - 1/2$ and $\text{AUK}_3 - 1/2$ measure negative dependence. Practical aspects of measurements for positive or negative dependence can be found, e.g., in Gargouri-Ellouze & Bargaoui (2009) and Eslamian (2014).

**The bounds.** The well-known Fréchet-Hoeffding result regarding copula bounds implies the following proposition.

**Proposition 2.** *For any continuous random variables $X$ and $Y$ distributed according to a bivariate distribution function $H(x,y)$, the measurements $\text{AUK}_i$, $i = 0, 1, 2, 3$, satisfy $1/4 \leq \text{AUK}_i \leq 1$, where the case with $H(x,y) = \Pr(X < x, X < y)$ provides $\text{AUK} = \text{AUK}_0$, $\text{AUK}_3$ to reach the lower bound $1/4$ and $\text{AUK}_1 = \text{AUK}_2 = 1$, the upper bound; whereas $H(x,y) = \Pr(X < x, -X < y)$ provides $\text{AUK}_0 = \text{AUK}_3 = 1$ and $\text{AUK}_1 = \text{AUK}_2 = 1/4$.*

**The AUK measure is countably sub-additive.** To formulate the next property we define the following AUKs based on the random vectors $(Z_1, W_1)$ and $(Z_2, W_2)$:

$$\text{AUK}_{i,0} = \Pr\{H_{i,0}(Z_i, W_i) < VU\}, \ i = 1, 2,$$

$$\text{AUK}_{12,0} = \Pr\{H_{12,0}(Z_1 + Z_2, W_1 + W_2) < VU\},$$

where the joint distribution functions $H_{i,0}(z,w) = \Pr(Z_i < z, W_i < w)$, $H_{12,0}(z,w) = \Pr(Z_1 + Z_2 < z, W_1 + W_2 < w)$ and $V, U$ are the independent and uniformly $[0,1]$ distributed random variables. The AUK measure is countably sub-additive, since we have the next proposition.

**Proposition 3.** *Let the random vectors $(Z_1, W_1)$, $(Z_2, W_2)$ be independent. Then, $\text{AUK}_{12,0} \leq \text{AUK}_{1,0} + \text{AUK}_{2,0}$.*

This statement also holds true, when the other components of the vector $\boldsymbol{D}$ are considered.

**Rényi's axioms.** According to Schweizer & Wolff (1981), we can present Rényi's conditions regarding a measure of dependence $R(X, Y)$ for two continuously distributed variables $X$ and $Y$ in the following form.

$R_1$: $R(X, Y)$ is defined for any $X$ and $Y$.

$R_2$: $R(X, Y) = R(Y, X)$.

$R_3$: $0 \leq R(X, Y) \leq 1$.

$R_4$: $R(X, Y) = 0$ if and only if $X$ and $Y$ are independent.

$R_5$: $R(X, Y) = 1$ if and only if each of $X$, $Y$ is a.s. a strictly monotone function of the other.

$R_6$: If $f$ and $g$ are strictly monotone a.s. on Range $X$ and Range $Y$, respectively, then $R\{f(X), g(Y)\} = R(X, Y)$.

$R_7$: If the joint distribution of $X$ and $Y$ is bivariate normal, with correlation coefficient $\rho$, then $R(X, Y)$ is a strictly increasing function of $|\rho|$.

$R_8$: If $(X, Y)$ and $(X_n, Y_n)$, $n = 1, 2, \ldots$, are pairs of random variables with joint distributions $H$ and $H_n$, respectively, and if the sequence $\{H_n\}$ converges weakly to $H$, then $\lim_{n \to \infty} R(X_n, Y_n) = R(X, Y)$.

It is clear that the proposed index $I_{\text{AUK}}$ satisfies requirements $R_1$, $R_2$ and $R_7$. Note that Proposition 2 offers bounds for $\text{AUK}_i, i = 0, \ldots, 3$. This proposition was employed to derive the measure $I_{\text{AUK}} \geq 0$ under the restriction that $I_{\text{AUK}} = 1$, when $(X, Y) = (X, X)$ or $(X, Y) = (X, -X)$.

In order to evaluate $I_{\text{AUK}}$ with respect to $R_4$, we state the assumption:

$C_1$: One can detect a subscript $i_0 \in [0, 1, 2, 3]$ such that the function $K_{i_0}(t) - \{t - t\log(t)\} \geq 0$ or $\leq 0$, for all $0 \leq t \leq 1$. That is, we can find at least one $K$-plot that does not cross the diagonal line.

Then the structure of the vector $\boldsymbol{D}$ used in the definition of $I_{\text{AUK}}$ plays a key role to prove the following proposition.

**Proposition 4.** *The measure $I_{\text{AUK}}$ satisfies $R_6$ and $R_8$. Moreover, under condition $C_1$, $I_{\text{AUK}}$ provides property $R_4$.*

*Remarks.* In general, there are cases of bivariate distribution functions that do not satisfy assumption $C_1$. For example, one can consider bivariate random vectors $(X, Y)'s$ that are uniformly distributed on the circumference of the circle $x^2 + y^2 = 1$. In this situation, we have $(X = \cos(\theta), Y = \sin(\theta))$, $\theta \sim Unif[0, 2\pi]$, and then, for all $i = 0, 1, 2, 3$,

$$K_i(t) = \begin{cases} t + 1/4, & t \in [0, 1/4), \\ 1, & t \in [1/4, 1]. \end{cases}$$

Note that $\overline{I}_{\text{AUK}}$ is a bijective increasing function of $I_{\text{AUK}}$; thus $\overline{I}_{\text{AUK}}$ provides the $I_{\text{AUK}}$-related results shown above.

The broad class of bivariate distribution functions can offer explicit or implicit forms of the function $K(t)$. We refer the reader to Genest & Rivest (1993, 2001); Nelsen et al. (2003) for corresponding examples. Genest & Rivest (1993) derived a general formula for computing the distribution function $K(t)$. Thus, for a broad class of distribution functions $H(x, y)$ one can calculate AUK in an analytical manner. In these cases, it is not difficult to show that using the remark presented in the last paragraph of Section 1 of Genest & Rivest (1993) we can also obtain analytically values of $\text{AUK}_i$, $i = 0, 1, 2, 3$. Although, in general, values of the vector $\boldsymbol{D} = (\text{AUK}_0, \ldots, \text{AUK}_3)^T$ can be easily computed numerically using, e.g., Proposition 1 (see the Web Supplementary Materials), we present the following example, when $\boldsymbol{D}$ can be expressed in an explicit form.

*Example.* For the Farlie-Gumbel-Morgenstern bivariate distribution function (1), by virtue of Proposition 1 we have

$$\text{AUK} = \text{E}[1 - C(V, U) + C(V, U) \log\{C(V, U)\}],$$

where $C(v, u) = vu\{1 + \gamma(1 - v)(1 - u)\}$, and the random variables $V$ and $U$ are uniformly $[0, 1]$ distributed. The joint distribution function $C(v, u)$ has the density function $(\partial^2 C(v, u))/(\partial v \partial u) = 1 + \gamma(1 - 2v)(1 - 2u)$. Therefore, the expression of the AUK as a function of $\gamma$ is

$$\text{AUK}(\gamma) = \int_0^1 \int_0^1 [1 - C(v, u) + C(v, u) \log\{C(v, u)\}]\{1 + \gamma(1 - 2v)(1 - 2u)\}\mathrm{d}v\mathrm{d}u.$$

This yields,

$$\begin{aligned}
\text{AUK}(\gamma) &= Re\big[\frac{1}{108}\gamma\{-29 + 6\text{Log}(1+\gamma)\} + \frac{1}{18\gamma}\{-13 + 3\pi^2 - 9\text{Log}(1+\gamma) \\
&\quad -18\text{Log}(-\gamma)\text{Log}(1+\gamma) - 18\text{Li}_2(1+\gamma)\} + \frac{1}{36\gamma^2}\{-\pi^2 + 20\text{Log}(1+\gamma) \\
&\quad +6\text{Log}(-\gamma)\text{Log}(1+\gamma) + 6\text{Li}_2(1+\gamma)\} \\
&\quad +\frac{1}{72}\{167 - 6\pi^2 - 72\text{Log}(1+\gamma) + 36\text{Log}(-\gamma)\text{Log}(1+\gamma) + 36\text{Li}_2(1+\gamma)\}\big]
\end{aligned}$$

with $\text{AUK}(0) = 1/2$, where the notation $\text{Li}_2(z)$ defines the dilogarithm function, $\text{Log}(z) = \log(z)$ if $z > 0$ and $\text{Log}(z) = \pi\imath + \log(-z)$ if $z < 0$, with $\imath$ to be the imaginary unit ($\imath^2 = -1$). A simple exercise in the use of the Wolfram Mathematica system of computer programs (Wolfram, 1999) for performing mathematical symbolic operations can confirm this expression and compute values of $\text{AUK}(\gamma)$ (see the Web Supplementary Materials for details). One can also express the AUK in the form

$$\text{AUK}(\gamma) = \text{E}([1 - C(V_1, U_1) + C(V_1, U_1)\log\{C(V_1, U_1)\}]\{1 + \gamma(1 - 2V_1)(1 - 2U_1)\}),$$

where the random variables $V_1$ and $U_1$ are independent and uniformly $[0, 1]$ distributed. It is convenient using the $R$ (R Development Core Team, 2008) function *"integrate"* or a Monte Carlo method to evaluate this AUK numerically. Similarly with the aforementioned explanations, one can show that $\text{AUK}_1(\gamma) = \text{AUK}_2(\gamma) = \text{AUK}(-\gamma)$ and $\text{AUK}_3(\gamma) = \text{AUK}_0(\gamma) = \text{AUK}(\gamma)$.

## 5. Simulation study

The goal of this section is two-fold. We first seek to provide empirical verification of the work presented above, and secondly to compare our proposed index with commonly used in practice indices, such as Pearson correlation and Kendall's $\tau$. An additional comparison we make is with the celebrated MIC index of dependence. To compute MIC we use the package "minerva" in $R$ (R Development Core Team, 2008).

**Normal Case:** We simulate $N = 1000$ random samples of size $n$ from the bivariate normal $N_2\left(\binom{0}{0}, \binom{1\ \rho}{\rho\ 1}\right)$ distribution, for various values of $n$ and $\rho$. Based on these samples we compute the averages (and the Monte Carlo standard deviations) of $I_{\text{AUK}}$, $\overline{I}_{\text{AUK}}$, MIC and $|r|$, where $r$ is the sample Pearson correlation. In the considered scenario, it is reasonable to assume that $|r|$, a parametric maximum likelihood based statistic, can demonstrate a standard in estimating the dependence between $X$ and $Y$. This design of Monte Carlo experiments can be used to judge the performance of the dependence measures, corresponding to values of $\rho$. The results are presented in Table 1. Observe that, as $n$ becomes large, $|r|$ and $\overline{I}_{\text{AUK}}$ are close to the true value of $|\rho|$, while the behavior of MIC is somewhat different. Specifically, MIC exhibits high bias for sample sizes up to 1000 when $\rho = 0$, while for $n = 5000$ the bias drops to 7.9%. Similarly, high biases are observed for MIC for all values of $|\rho| < 1$ and all sample sizes. On the other hand, $\overline{I}_{\text{AUK}}$ closely follows the true value of $|\rho|$, where higher bias is observed for small sample sizes (less than 100) and values of $|\rho| \leq 0.1$. Therefore, $\overline{I}_{\text{AUK}}$ in the case of normal distribution approximately behaves as $|\rho|$, when $n > 100$. Notice also that $|r|$ is a measure of general dependence for the case of bivariate normal distributions. On the other hand, $\overline{I}_{\text{AUK}}$ is a measure of dependence for general continuous marginal distributions, not only necessarily for bivariate normal distributions. Therefore, $\overline{I}_{\text{AUK}}$ is a general index of dependence and it presents linearly the dependence structure of a bivariate normal distribution function.

It is interesting to note that $I_{\text{AUK}} \cong 0$ and $\overline{I}_{\text{AUK}} \cong 0$ when $\rho = 0$ for $n \geq 50$. (In these cases, the symbol '$\cong$' is used, taking into account corresponding values of $|r| - |\rho|$, where $|r|$ is the maximum likelihood estimator of $|\rho|$.) On the other hand, MIC is clearly not close to 0 when $\rho = 0$ and it is sensitive to the size of the sample. This can be easily seen from Table 1; when the sample size is 1000 MIC equals 0.133 exhibiting a 13.3% bias (when $\rho = 0$), while when $n = 5000$ the bias is only 7.9% (MIC = 0.079), as we have discussed above. It seems that $I_{\text{AUK}}$ is somewhat better than $\overline{I}_{\text{AUK}}$ and MIC with respect to the distances $|I_{\text{AUK}} - |\rho||$, $|\overline{I}_{\text{AUK}} - |\rho||$ and $|\text{MIC} - |\rho||$, when $\rho < 0.2$ and $n \leq 100$.

The results shown above regarding the scenarios with $\rho = 0$ (observed vectors consist of independent components, however $\text{MIC} > \max(I_{\text{AUK}}, \overline{I}_{\text{AUK}}))$ raise a concern regarding the conclusion: "if $\text{MIC} > I_{\text{AUK}}$ (or $\text{MIC} > \overline{I}_{\text{AUK}}$), then MIC outperforms $I_{\text{AUK}}$ (or $\overline{I}_{\text{AUK}}$, respectively)".

**Non-Normal Cases:** Table 2 shows the average values (and the Monte Carlo standard deviations) of $I_{\text{AUK}}$, $\overline{I}_{\text{AUK}}$, MIC and $|\tau|$, where $\tau$ is Kendall's $\tau$, a well-accepted nonparametric index. The considered measures were obtained using random samples $(X_1, Y_1), \ldots, (X_n, Y_n)$ from various bivariate distributions. The number of Monte Carlo repetitions is $N = 5000$. We refer the reader to Johnson (2013) for details regarding the definitions of the bivariate distributions that are listed in Table 2. Johnson (2013) provided simple algorithms related to the appropriate calculations for generating samples from the considered bivariate distributions. In order to obtain the Monte Carlo results based on random variables from a bivariate $t_5$ distribution we use the package "mvtnorm" in $R$ (R Development Core Team, 2008), generating data from the bivariate student distribution with 5 degree of freedom and variance-covariance matrix $\begin{pmatrix} 1 & 1 \\ 1 & 4 \end{pmatrix}$. We also

Table 1: The Monte Carlo expectations (and the Monte Carlo standard deviations) of $|r|$, $I_{\text{AUK}}$, $\bar{I}_{\text{AUK}}$ and MIC; the number of Monte Carlo repetitions is $N = 1000$. Random samples of size $n$ are drawn from a bivariate normal distribution with correlation coefficient $\rho$, for various values of $n$ and $\rho$.

|  |  | \|$\rho$\| | | | | | |
|---|---|---|---|---|---|---|---|
|  |  | 0 | 0.1 | 0.2 | 0.3 | 0.4 | 0.5 |
| $n = 20$ | $\|r\|$ | 0.193(0.141) | 0.198(0.144) | 0.241(0.160) | 0.314(0.176) | 0.388(0.184) | 0.489(0.172) |
| | $I_{\text{AUK}}$ | 0.157(0.067) | 0.169(0.074) | 0.198(0.085) | 0.233(0.095) | 0.270(0.101) | 0.326(0.102) |
| | $\bar{I}_{\text{AUK}}$ | 0.313(0.122) | 0.334(0.134) | 0.386(0.150) | 0.447(0.160) | 0.508(0.164) | 0.595(0.151) |
| | MIC | 0.361(0.104) | 0.364(0.107) | 0.370(0.107) | 0.386(0.120) | 0.415(0.127) | 0.465(0.151) |
| $n = 30$ | $\|r\|$ | 0.145(0.110) | 0.167(0.126) | 0.215(0.142) | 0.310(0.156) | 0.400(0.152) | 0.490(0.144) |
| | $I_{\text{AUK}}$ | 0.112(0.051) | 0.131(0.064) | 0.159(0.074) | 0.204(0.080) | 0.250(0.084) | 0.299(0.085) |
| | $\bar{I}_{\text{AUK}}$ | 0.228(0.099) | 0.264(0.122) | 0.317(0.137) | 0.399(0.143) | 0.479(0.142) | 0.558(0.134) |
| | MIC | 0.261(0.075) | 0.264(0.074) | 0.279(0.080) | 0.298(0.086) | 0.331(0.099) | 0.375(0.115) |
| $n = 50$ | $\|r\|$ | 0.116(0.087) | 0.144(0.102) | 0.202(0.117) | 0.294(0.122) | 0.395(0.123) | 0.497(0.109) |
| | $I_{\text{AUK}}$ | 0.080(0.037) | 0.101(0.051) | 0.132(0.058) | 0.175(0.063) | 0.228(0.068) | 0.283(0.066) |
| | $\bar{I}_{\text{AUK}}$ | 0.164(0.075) | 0.205(0.101) | 0.266(0.112) | 0.349(0.117) | 0.443(0.119) | 0.536(0.108) |
| | MIC | 0.297(0.057) | 0.303(0.058) | 0.313(0.060) | 0.336(0.068) | 0.369(0.082) | 0.412(0.089) |
| $n = 100$ | $\|r\|$ | 0.080(0.060) | 0.114(0.078) | 0.203(0.096) | 0.294(0.094) | 0.399(0.085) | 0.500(0.077) |
| | $I_{\text{AUK}}$ | 0.050(0.026) | 0.071(0.038) | 0.115(0.048) | 0.159(0.049) | 0.214(0.047) | 0.271(0.045) |
| | $\bar{I}_{\text{AUK}}$ | 0.104(0.054) | 0.146(0.077) | 0.233(0.094) | 0.319(0.094) | 0.422(0.085) | 0.520(0.076) |
| | MIC | 0.236(0.034) | 0.243(0.036) | 0.255(0.042) | 0.274(0.047) | 0.313(0.053) | 0.359(0.065) |
| $n = 200$ | $\|r\|$ | 0.053(0.039) | 0.103(0.064) | 0.198(0.067) | 0.299(0.066) | 0.399(0.058) | 0.499(0.054) |
| | $I_{\text{AUK}}$ | 0.033(0.018) | 0.059(0.031) | 0.103(0.034) | 0.154(0.035) | 0.207(0.032) | 0.263(0.032) |
| | $\bar{I}_{\text{AUK}}$ | 0.068(0.036) | 0.121(0.064) | 0.212(0.068) | 0.310(0.067) | 0.410(0.059) | 0.508(0.055) |
| | MIC | 0.221(0.024) | 0.225(0.024) | 0.237(0.028) | 0.263(0.032) | 0.295(0.036) | 0.342(0.041) |
| $n = 300$ | $\|r\|$ | 0.047(0.036) | 0.104(0.054) | 0.199(0.056) | 0.299(0.054) | 0.400(0.049) | 0.497(0.043) |
| | $I_{\text{AUK}}$ | 0.028(0.016) | 0.056(0.026) | 0.102(0.028) | 0.152(0.029) | 0.205(0.027) | 0.260(0.026) |
| | $\bar{I}_{\text{AUK}}$ | 0.057(0.034) | 0.116(0.054) | 0.209(0.056) | 0.307(0.056) | 0.408(0.050) | 0.504(0.044) |
| | MIC | 0.192(0.018) | 0.198(0.018) | 0.210(0.021) | 0.235(0.026) | 0.271(0.031) | 0.315(0.034) |
| $n = 500$ | $\|r\|$ | 0.038(0.028) | 0.100(0.044) | 0.200(0.041) | 0.301(0.040) | 0.400(0.038) | 0.500(0.033) |
| | $I_{\text{AUK}}$ | 0.022(0.012) | 0.052(0.021) | 0.101(0.021) | 0.151(0.021) | 0.204(0.021) | 0.259(0.020) |
| | $\bar{I}_{\text{AUK}}$ | 0.046(0.026) | 0.108(0.044) | 0.206(0.042) | 0.305(0.040) | 0.405(0.039) | 0.503(0.034) |
| | MIC | 0.163(0.012) | 0.167(0.013) | 0.181(0.016) | 0.205(0.018) | 0.240(0.022) | 0.287(0.026) |
| $n = 1000$ | $\|r\|$ | 0.025(0.020) | 0.099(0.031) | 0.198(0.030) | 0.300(0.029) | 0.401(0.027) | 0.499(0.024) |
| | $I_{\text{AUK}}$ | 0.015(0.009) | 0.050(0.015) | 0.098(0.015) | 0.149(0.016) | 0.202(0.015) | 0.258(0.015) |
| | $\bar{I}_{\text{AUK}}$ | 0.030(0.019) | 0.103(0.032) | 0.201(0.031) | 0.302(0.030) | 0.403(0.028) | 0.501(0.025) |
| | MIC | 0.133(0.008) | 0.138(0.008) | 0.152(0.010) | 0.176(0.013) | 0.211(0.015) | 0.257(0.018) |
| $n = 5000$ | $\|r\|$ | 0.011(0.008) | 0.100(0.014) | 0.200(0.014) | 0.300(0.013) | 0.400(0.012) | 0.500(0.010) |
| | $I_{\text{AUK}}$ | 0.006(0.004) | 0.048(0.007) | 0.098(0.007) | 0.148(0.007) | 0.201(0.006) | 0.257(0.006) |
| | $\bar{I}_{\text{AUK}}$ | 0.013(0.008) | 0.100(0.014) | 0.201(0.014) | 0.300(0.013) | 0.400(0.012) | 0.501(0.011) |
| | MIC | 0.079(0.003) | 0.084(0.003) | 0.098(0.004) | 0.121(0.005) | 0.157(0.006) | 0.204(0.007) |

Table 1: (continued)

|  |  | $|\rho|$ | | | | |
|---|---|---|---|---|---|---|
|  |  | 0.6 | 0.7 | 0.8 | 0.9 | 1 |
| $n=20$ | $|r|$ | 0.594(0.151) | 0.697(0.121) | 0.793(0.090) | 0.893(0.051) | 1(0) |
|  | $I_{\text{AUK}}$ | 0.383(0.104) | 0.453(0.095) | 0.530(0.096) | 0.641(0.082) | 1(0) |
|  | $\bar{I}_{\text{AUK}}$ | 0.674(0.139) | 0.762(0.108) | 0.836(0.087) | 0.916(0.047) | 1(0) |
|  | MIC | 0.522(0.153) | 0.607(0.163) | 0.683(0.160) | 0.811(0.128) | 1(0) |
| $n=30$ | $|r|$ | 0.593(0.126) | 0.695(0.093) | 0.793(0.072) | 0.895(0.038) | 1(0) |
|  | $I_{\text{AUK}}$ | 0.358(0.086) | 0.429(0.076) | 0.510(0.073) | 0.625(0.062) | 1(0) |
|  | $\bar{I}_{\text{AUK}}$ | 0.645(0.121) | 0.739(0.091) | 0.825(0.068) | 0.911(0.037) | 1(0) |
|  | MIC | 0.432(0.127) | 0.517(0.131) | 0.615(0.137) | 0.762(0.114) | 1(0) |
| $n=50$ | $|r|$ | 0.602(0.092) | 0.694(0.077) | 0.798(0.054) | 0.900(0.028) | 1(0) |
|  | $I_{\text{AUK}}$ | 0.347(0.064) | 0.409(0.060) | 0.499(0.055) | 0.616(0.046) | 1(0) |
|  | $\bar{I}_{\text{AUK}}$ | 0.636(0.092) | 0.720(0.076) | 0.819(0.053) | 0.909(0.029) | 1(0) |
|  | MIC | 0.482(0.102) | 0.553(0.106) | 0.660(0.107) | 0.807(0.092) | 1(0) |
| $n=100$ | $|r|$ | 0.601(0.066) | 0.701(0.051) | 0.799(0.037) | 0.899(0.019) | 1(0) |
|  | $I_{\text{AUK}}$ | 0.333(0.045) | 0.402(0.041) | 0.486(0.038) | 0.605(0.031) | 1(0) |
|  | $\bar{I}_{\text{AUK}}$ | 0.618(0.068) | 0.713(0.052) | 0.809(0.039) | 0.904(0.019) | 1(0) |
|  | MIC | 0.423(0.070) | 0.502(0.072) | 0.605(0.075) | 0.757(0.071) | 1(0) |
| $n=200$ | $|r|$ | 0.599(0.045) | 0.700(0.037) | 0.799(0.026) | 0.899(0.014) | 1(0) |
|  | $I_{\text{AUK}}$ | 0.325(0.031) | 0.396(0.030) | 0.481(0.026) | 0.599(0.022) | 1(0) |
|  | $\bar{I}_{\text{AUK}}$ | 0.607(0.047) | 0.707(0.038) | 0.805(0.027) | 0.901(0.014) | 1(0) |
|  | MIC | 0.404(0.046) | 0.486(0.053) | 0.593(0.053) | 0.747(0.050) | 1(0) |
| $n=300$ | $|r|$ | 0.600(0.039) | 0.699(0.030) | 0.799(0.021) | 0.900(0.011) | 1(0) |
|  | $I_{\text{AUK}}$ | 0.324(0.026) | 0.393(0.024) | 0.478(0.021) | 0.598(0.018) | 1(0) |
|  | $\bar{I}_{\text{AUK}}$ | 0.606(0.040) | 0.704(0.031) | 0.802(0.022) | 0.900(0.011) | 1(0) |
|  | MIC | 0.380(0.038) | 0.459(0.042) | 0.567(0.043) | 0.724(0.042) | 1(0) |
| $n=500$ | $|r|$ | 0.600(0.030) | 0.700(0.023) | 0.799(0.016) | 0.900(0.009) | 1(0) |
|  | $I_{\text{AUK}}$ | 0.321(0.020) | 0.392(0.018) | 0.477(0.016) | 0.596(0.014) | 1(0) |
|  | $\bar{I}_{\text{AUK}}$ | 0.603(0.030) | 0.702(0.023) | 0.801(0.017) | 0.900(0.009) | 1(0) |
|  | MIC | 0.349(0.030) | 0.430(0.031) | 0.540(0.034) | 0.699(0.034) | 1(0) |
| $n=1000$ | $|r|$ | 0.600(0.020) | 0.700(0.016) | 0.800(0.012) | 0.900(0.006) | 1(0) |
|  | $I_{\text{AUK}}$ | 0.320(0.014) | 0.391(0.013) | 0.476(0.012) | 0.595(0.009) | 1(0) |
|  | $\bar{I}_{\text{AUK}}$ | 0.601(0.021) | 0.701(0.017) | 0.801(0.012) | 0.899(0.006) | 1(0) |
|  | MIC | 0.320(0.020) | 0.403(0.022) | 0.514(0.023) | 0.674(0.024) | 1(0) |
| $n=5000$ | $|r|$ | 0.600(0.009) | 0.700(0.007) | 0.800(0.005) | 0.900(0.003) | 1(0) |
|  | $I_{\text{AUK}}$ | 0.319(0.006) | 0.390(0.006) | 0.476(0.005) | 0.595(0.004) | 1(0) |
|  | $\bar{I}_{\text{AUK}}$ | 0.600(0.009) | 0.700(0.007) | 0.800(0.005) | 0.899(0.003) | 1(0) |
|  | MIC | 0.267(0.009) | 0.350(0.010) | 0.462(0.010) | 0.625(0.010) | 1(0) |

evaluate the case with $(X, Y = \varepsilon/X^2)$, where $X$ and $\varepsilon$ are independent $N(5, 1)$ distributed random variables. This case can illustrate a practical issue related to measurement error problems (e.g., Vexler et al., 2014), when dependence measures based on $Cov(X, Y)$ are improper.

Consider Tables 2a and 2b. In order to generate data points $(X, Y)$'s from a bivariate (BV) Morgenstern type distribution, we use the following scheme. Let $X, U$ be independent $Unif[0, 1]$ distributed. Define $Z = \alpha(2X - 1) - 1$, $W = 1 - 2\alpha(2X - 1) + \alpha^2(2X - 1)^2 + 4\alpha U(2X - 1)$. Then $Y = 2U/(W^{1/2} - Z)$. To facilitate discussion of the results presented in Tables 2a and 2b, Figure 4 presents illustrations of graphs of a single sample $\{(X_1, Y_1), ..., (X_n, Y_n)\}$ of size $n = 1000$ drawn from the BV Morgenstern distribution with parameters $\alpha = 0.5$ and $\alpha = 5$, respectively. Figure 5 shows the corresponding multi-panel K-plots related to the BV Morgenstern distributions with $\alpha = 0.5$ and $\alpha = 5$. These graphs depict an increase of the dependence between $X$ and $Y$ when the parameter $\alpha$ is changed from 0.5 to 5 (see Johnson (2013) for more details, in this context). The index $|\tau|$ shows a relative increase of about $(0.6 - 0.12)/0.12 = 4$, the index $I_{\text{AUK}}$ shows a relative increase of about $(0.5 - 0.1)/0.1 = 4$, the index $\overline{I}_{\text{AUK}}$ shows a relative increase of about $(0.8 - 0.2)/0.2 = 3$, whereas the index MIC provides only a relative increase of about $(0.85 - 0.25)/0.25 = 2.4$. (These calculations are provided approximately taking into account values of the indexes with respect to $n = 25, \ldots, 200$.) Note that, in Figure 5, the K-plots indicate that variables $(X, Y)$'s $\sim$ BV-Morgenstern distribution with $\alpha = 5$ are "more positively" dependent than variables $(X, Y)$'s $\sim$ BV-Morgenstern distribution with $\alpha = 0.5$.

Consider Tables 2c and 2d. In this case, to generate data points $(X, Y)$'s from the BV Plackett distribution, we use the following scheme. Let $X$ and $U$ be independent $Unif[0, 1]$ distributed. Define $W_1 = U(1-U)$, $W_2 = \psi + W_1(\psi - 1)^2$, $W_3 = 2W_1(\psi^2 X + 1 - X) + \psi(1 - 2W_1)$, $W_4 = \psi[\psi + 4(1 - \psi)^2 X(1 - X)W_1]$. Then $Y = W_2[W_3 - (1 - 2U)W_4^{1/2}]/2$. The parameter $\psi$ of the BV Plackett distribution can characterize the dependence between $X$ and $Y$ (see Johnson (2013) for more details, in this context). Regarding Tables 2c and 2d, we can observe that the index $|\tau|$ shows a relative increase of about $(0.18 - 0.09)/0.09 = 1$, the index $I_{\text{AUK}}$ shows a relative increase of about $(0.17 - 0.07)/0.07 \simeq 1.4$, the index $\overline{I}_{\text{AUK}}$ shows a relative increase of about $(0.3 - 0.13)/0.13 \simeq 1.3$, whereas the index MIC provides only a relative increase of about $(0.3 - 0.25)/0.25 = 0.2$.

Perhaps, the results shown above can be employed to state that, in the considered scenarios, the proposed indexes detect the change in the dependence structure of observations more accurately and sensitively when compared with the MIC.

Tables 2g and 2h clearly indicate that the proposed indexes are able to recognize dependence for all cases studied providing accurate measurements of the degree of dependence in samples with relatively small and large sizes coming from a variety of distributions.

Note that, the goal of this article is not directly related to outperforming the classical measures of dependence or MIC. The proposed approach provides a meaningful graphical method for understanding bivariate dependence that is demonstrated in Section 6.

Table 2: The Monte Carlo expectations (and the Monte Carlo standard deviations) of $|\tau|$ ($\tau$ is the Kendall's $\tau$), $I_{\text{AUK}}$, $\bar{I}_{\text{AUK}}$ and MIC; the number of Monte Carlo repetitions is $N = 5000$. Random samples of size $n$ are drawn from different distributions, for various values of $n$.

(a) The BV Morgenstern type distribution, $\alpha = 0.5$.

|  | $n$ | | | |
|---|---|---|---|---|
|  | 25 | 50 | 100 | 200 |
| $|\tau|$ | .148(.102) | .123(.079) | .114(.061) | .112(.046) |
| $I_{\text{AUK}}$ | .167(.075) | .121(.057) | .100(.044) | .090(.034) |
| $\bar{I}_{\text{AUK}}$ | .331(.138) | .246(.111) | .205(.087) | .185(.068) |
| MIC | .315(.091) | .311(.062) | .253(.040) | .236(.027) |

(b) The BV Morgenstern type distribution, $\alpha = 5$.

|  | $n$ | | | |
|---|---|---|---|---|
|  | 25 | 50 | 100 | 200 |
| $|\tau|$ | .621(.072) | .622(.046) | .623(.030) | .622(.021) |
| $I_{\text{AUK}}$ | .532(.068) | .510(.044) | .499(.030) | .493(.020) |
| $\bar{I}_{\text{AUK}}$ | .845(.059) | .830(.040) | .822(.028) | .817(.020) |
| MIC | .836(.103) | .873(.083) | .854(.063) | .857(.045) |

(c) The BV Plackett distribution, $\psi = 1.25$.

|  | $n$ | | | |
|---|---|---|---|---|
|  | 25 | 50 | 100 | 200 |
| $|\tau|$ | .120(.090) | .088(.064) | .068(.049) | .058(.039) |
| $I_{\text{AUK}}$ | .142(.063) | .092(.047) | .066(.036) | .051(.028) |
| $\bar{I}_{\text{AUK}}$ | .286(.119) | .188(.093) | .136(.073) | .106(.058) |
| MIC | .303(.084) | .299(.056) | .240(.036) | .223(.024) |

(d) The BV Plackett distribution, $\psi = 2$.

|  | $n$ | | | |
|---|---|---|---|---|
|  | 25 | 50 | 100 | 200 |
| $|\tau|$ | .202(.121) | .189(.093) | .187(.068) | .187(.047) |
| $I_{\text{AUK}}$ | .220(.087) | .181(.067) | .164(.049) | .156(.035) |
| $\bar{I}_{\text{AUK}}$ | .425(.151) | .359(.125) | .328(.094) | .315(.067) |
| MIC | .338(.098) | .345(.070) | .292(.046) | .283(.032) |

(e) ALI HAQ's BV distribution with $a = 0.9$ and $p = 0.5$: Johnson (2013, p. 202).

|  | $n$ | | | |
|---|---|---|---|---|
|  | 25 | 50 | 100 | 200 |
| $|\tau|$ | .671(.092) | .671(.062) | .670(.043) | .671(.030) |
| $I_{\text{AUK}}$ | .600(.084) | .580(.057) | .569(.040) | .565(.027) |
| $\bar{I}_{\text{AUK}}$ | .892(.058) | .884(.041) | .879(.029) | .877(.020) |
| MIC | .756(.131) | .771(.100) | .729(.074) | .721(.053) |

(f) Gumbel's BV exponential distribution, $\theta = 0.9$: Johnson (2013, p. 197).

|  | $n$ | | | |
|---|---|---|---|---|
|  | 25 | 50 | 100 | 200 |
| $|\tau|$ | .164(.109) | .144(.085) | .139(.066) | .139(.047) |
| $I_{\text{AUK}}$ | .192(.076) | .149(.059) | .128(.045) | .119(.033) |
| $\bar{I}_{\text{AUK}}$ | .377(.136) | .299(.112) | .259(.088) | .243(.065) |
| MIC | .330(.099) | .327(.066) | .271(.045) | .256(.031) |

(g) BV-$t_5$ distribution.

|  | $n$ | | | |
|---|---|---|---|---|
|  | 25 | 50 | 100 | 200 |
| $|\tau|$ | .336(.132) | .333(.094) | .333(.065) | .333(.045) |
| $I_{\text{AUK}}$ | .314(.100) | .284(.073) | .270(.051) | .265(.036) |
| $\bar{I}_{\text{AUK}}$ | .577(.152) | .537(.119) | .519(.086) | .511(.061) |
| MIC | .410(.127) | .409(.090) | .350(.060) | .335(.042) |

(h) $(X, Y = \varepsilon/X^2)$, where $X$ and $\varepsilon$ are independent $N(5, 1)$ distributions.

|  | $n$ | | | |
|---|---|---|---|---|
|  | 25 | 50 | 100 | 200 |
| $|\tau|$ | .700(.078) | .700(.052) | .700(.035) | .700(.024) |
| $I_{\text{AUK}}$ | .516(.071) | .544(.049) | .562(.034) | .572(.023) |
| $\bar{I}_{\text{AUK}}$ | .831(.068) | .859(.041) | .874(.026) | .882(.017) |
| MIC | .783(.135) | .786(.097) | .741(.074) | .730(.052) |

## 6. Real data set example related to biomarkers in Myocardial Infarction (MI)

In order to illustrate the applicability of the proposed work we use real life data example that is based on a sample from a study that evaluates biomarkers related to the MI disease. The study was focused on the residents of Erie and Niagara counties, 35–79 years of age. The New

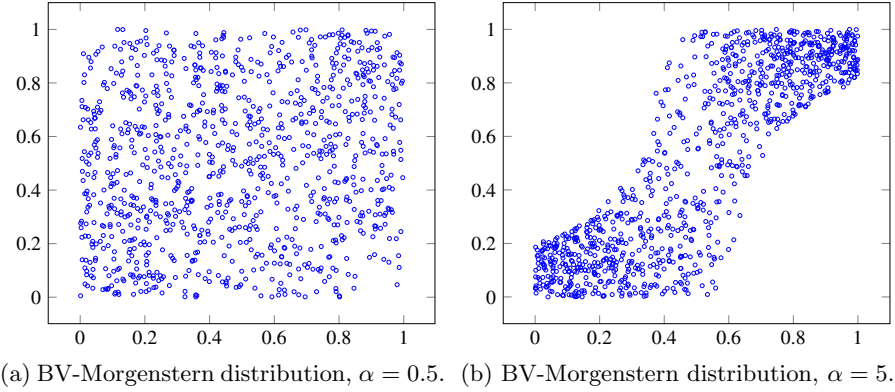

(a) BV-Morgenstern distribution, $\alpha = 0.5$. (b) BV-Morgenstern distribution, $\alpha = 5$.

Figure 4: Scatterplots: Random samples of size $n = 1000$ are drawn from the bivariate (BV) distributions.

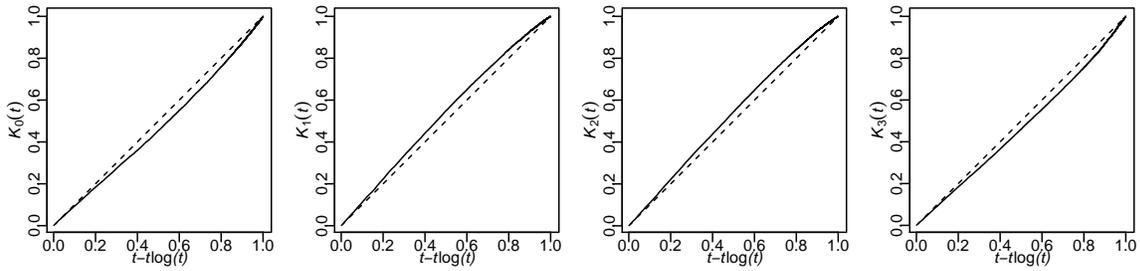

(a) BV-Morgenstern distribution, $\alpha = 0.5$.

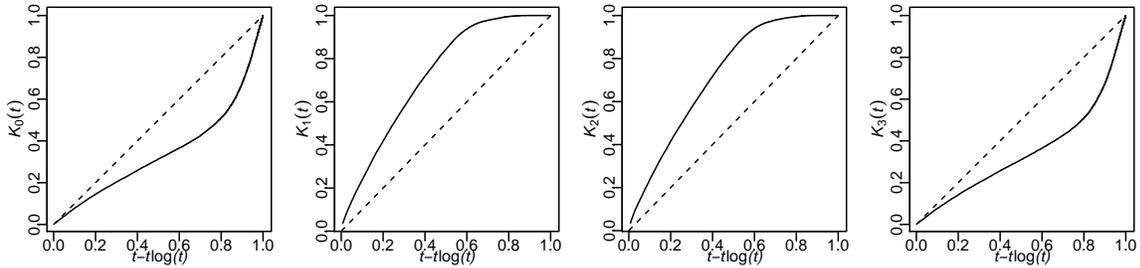

(b) BV-Morgenstern distribution, $\alpha = 5$.

Figure 5: The multi-panel K-plots related to the BV Morgenstern distributions with parameters $\alpha = 0.5$ and $\alpha = 5$, respectively.

York State department of Motor Vehicles drivers' license rolls were used as the sampling frame for adults between the age of 35 and 65 years, while the elderly sample (age 65–79) was randomly chosen from the Health Care Financing Administration database. We consider the biomarker "thiobarbituric acid-reactive substances" (TBARS). TBARS is commonly used to summarize the antioxidant status process of an individual in laboratory research (Armstrong, 2012), but its use as a discriminant factor between individuals with and without MI disease is still controversial (e.g., Schisterman et al., 2001). In the study investigating the discriminant ability of TBARS with regard to MI disease, dependencies between TBARS and other antioxidant status measures related to MI disease are evaluated. The literature has widely addressed concerns regarding assumptions for fitting various parametric distribution functions, including normal distributions, to actual TBARS' distributions. A number of antioxidants were examined from fresh blood samples at baseline, including TBARS and high-density lipoprotein (HDL) cholesterol. The HDL-cholesterol biomarker is often used as a good discriminant factor between individuals with and without MI disease (e.g., Schisterman et al., 2001).

In this study we measure dependence between TBARS and HDL-cholesterol, with a sample of 545 individuals with the MI disease (MI = 1) and a sample of 1495 individuals without the MI disease (MI = 0). Figure 6 depicts the scatterplots based on observed values of (TBARS,HDL) for MI = 0 and MI = 1 cases, respectively. Figure 7 presents the multi-panel K-plots of TBARS and HDL-cholesterol for both MI = 0 and MI = 1 cases. The corresponding K-plots were approximated using the empirical estimators of the probability functions $H_i$ and $\Pr\{H_i(X,Y) < t\}$, $i = 0, 1, 2, 3$, defined in Section 3.

The first panel of the K-plot indicates the type of dependence between the two biomarkers TBARS and HDL. Note that in both cases, i.e. MI = 0 and MI = 1, the first panel of the K-plot is slightly above the diagonal line, indicating that TBARS and HDL are negatively dependent. It is of interest to note that the Pearson correlation coefficient is $r = -0.11$ (MI = 0), value that agrees with the first Kendall plot (MI = 0), while $r = 0.04$ (MI = 1), value that disagrees with the negative dependence in the first Kendall plot for the case MI = 1.

We also note that despite the negative dependence indicated by the first two panels of the K-plots, the corresponding curves are very close to the diagonal line which indicates that the two variables are independent. In fact, the remaining panels of the K-plots corroborate the near independence of TBARS and HDL biomarkers.

To conduct nonparametric estimation of components of the AUKs based vector $\boldsymbol{D}$, the nonparametric technique introduced in Section 4 was executed. The respective bootstrap 95% Con-fidence Intervals for the measures were obtained via resampling 5000 times from the biomarkers' values. Table 3 presents our results.

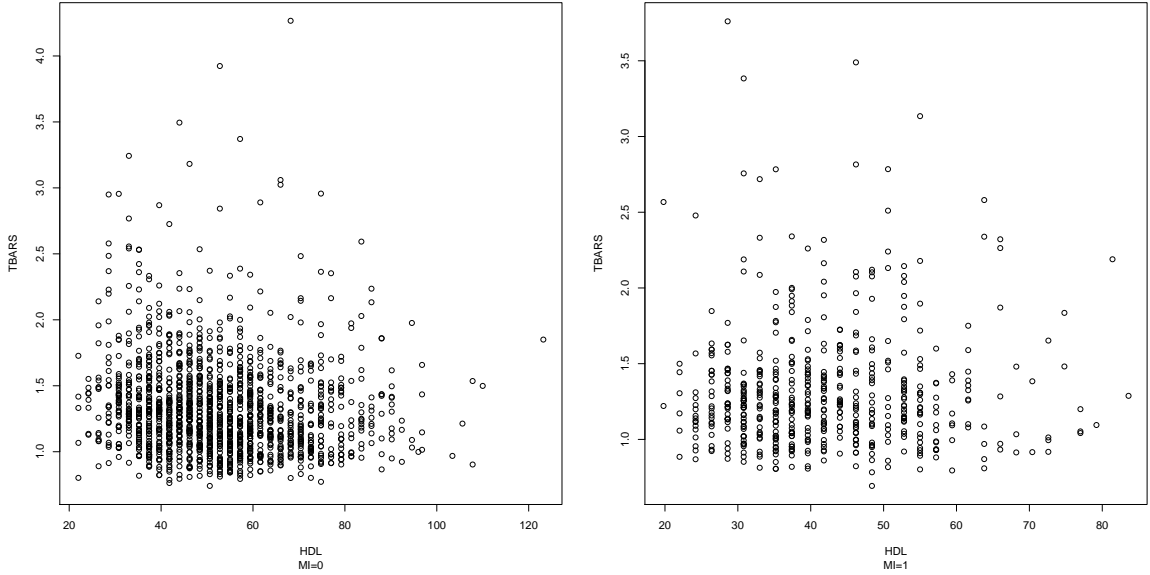

Figure 6: The scatterplots TBARS vs. HDL, for MI = 0 and MI = 1 cases.

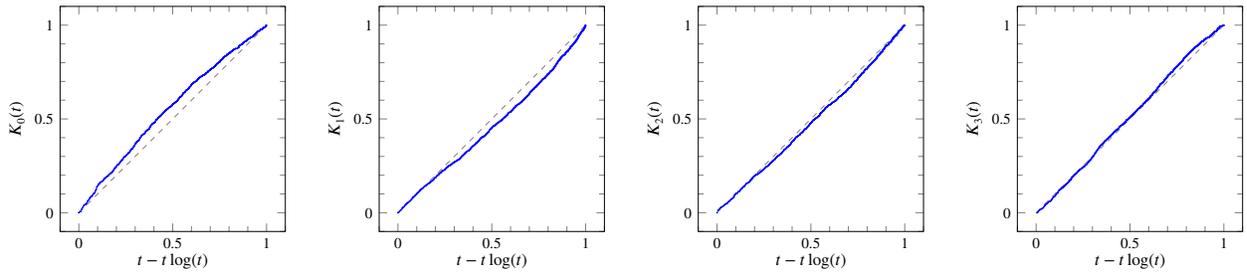

(a) The multi-panel K-plot for MI = 0 case.

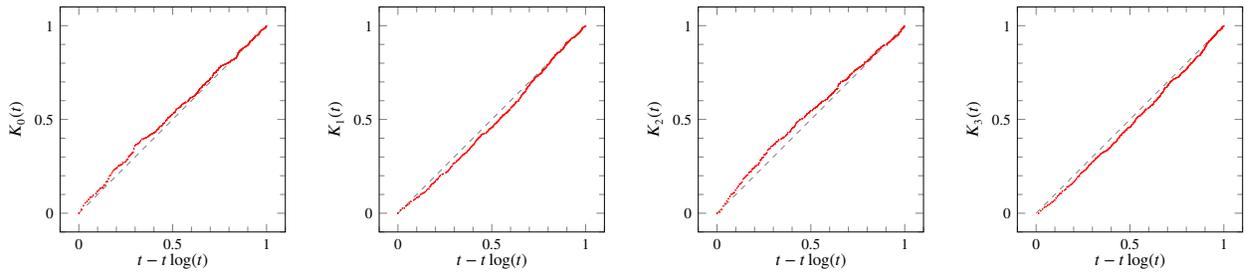

(b) The multi-panel K-plot for MI = 1 case.

Figure 7: The multi-panel K-plots of TBARS, HDL-cholesterol.

Table 3: Estimated Associations between TBARS and HDL-cholesterol.

|  | Estimator | 90% CI | 95% CI |
| --- | --- | --- | --- |
|  |  | MI = 0 |  |
| MIC | 0.1152 | [0.1698, 0.2096] | [0.1657, 0.2135] |
| $AUK_0$ | 0.5528 | [0.5443, 0.5623] | [0.5422, 0.5639] |
| $AUK_1$ | 0.4634 | [0.4556, 0.4720] | [0.4542, 0.4736] |
| $AUK_2$ | 0.4820 | [0.4737, 0.4911] | [0.4722, 0.4929] |
| $AUK_3$ | 0.5108 | [0.5010, 0.5193] | [0.4991, 0.5206] |
| $I_{AUK}$ | 0.0855 | [0.0676, 0.1041] | [0.0649, 0.1079] |
| $\overline{I}_{AUK}$ | 0.1760 | [0.1395, 0.2135] | [0.1339, 0.2211] |
|  |  | MI = 1 |  |
| MIC | 0.1381 | [0.1986, 0.2620] | [0.1934, 0.2684] |
| $AUK_0$ | 0.5245 | [0.5102, 0.5408] | [0.5078, 0.5439] |
| $AUK_1$ | 0.4791 | [0.4661, 0.4951] | [0.4635, 0.4972] |
| $AUK_2$ | 0.5313 | [0.5168, 0.5489] | [0.5138, 0.5514] |
| $AUK_3$ | 0.4752 | [0.4580, 0.4872] | [0.4556, 0.4904] |
| $I_{AUK}$ | 0.0649 | [0.0593, 0.0846] | [0.0571, 0.0880] |
| $\overline{I}_{AUK}$ | 0.1339 | [0.1225, 0.1741] | [0.1179, 0.1810] |

Notice that the 90% confidence intervals for $AUK_0$, $AUK_2$ and $AUK_3$ for MI = 0 and MI = 1 do not overlap. This is due to the fact that there are slight deviations of the Kendall curve from the diagonal line that are translated into non-overlapping confidence intervals for $AUK_0$, $AUK_2$ and $AUK_3$. Figure 7 clearly indicates that the Kendall curve $K_2(t)$ is slightly below the diagonal line when MI = 0, while $K_2(t)$ is above the diagonal line when MI = 1. Similar observations apply to $K_3(t)$ for MI = 0 and MI = 1. On the other hand, there is considerable overlap between the confidence intervals for $\overline{I}_{AUK}$ when MI = 0 and MI = 1 indicating that the same degree of dependence observed in Figure 7 holds for both conditions. The fact that the confidence intervals for $AUK_0$, $AUK_2$ and $AUK_3$ are non-overlapping is clearly confirmed for the cases of $AUK_2$ and $AUK_3$ by the 97.5% confidence intervals presented in Table 4. The confidence intervals for $AUK_j, j = 0, 1, 2, 3$, shown in Tables 3 and 4 provide significant arguments in favor of a negative correlation between values of the biomarkers TBARS and HDL, when MI = 0, 1. The corresponding 95% confidence intervals for the Pearson correlation coefficient are $[-0.1610, -0.0609]$ and $[-0.0582, 0.1197]$, for MI = 0 and MI = 1, respectively.

Table 4: The 97.5% CI estimators of AUKs.

|  | MI = 0 | MI = 1 |
| --- | --- | --- |
| $AUK_0$ | [0.5412, 0.5651] | [0.5049, 0.5469] |
| $AUK_1$ | [0.4530, 0.4750] | [0.4611, 0.5000] |
| $AUK_2$ | [0.4706, 0.4947] | [0.5108, 0.5544] |
| $AUK_3$ | [0.4975, 0.5224] | [0.4530, 0.4928] |

According to the conclusions made in Section 4, in the case MI = 0, the 90% confidence interval for $AUK_0$ depicted in Table 3 shows evidence of a "more negative" dependence between TBARS and HDL, in comparison with the case MI = 1.

Thus, the proposed **D**-vector based approach detects a difference in dependency structures related to joint distributions of (TBARS, HDL-cholesterol) with respect to MI = 0 and MI = 1. The index MIC does not provide significant results to discriminate values of (TBARS, HDL-cholesterol) between "disease" and "non-disease" populations.

## 7. Concluding remarks and discussion

In this paper we employ the ROC/AUC methodology in order to propose the multi-panel Kendall plot as an extension of the K-plot to measure dependence between two continuous random variables. This extension was necessitated by the inability of the K-plot to sufficiently represent a wide range of two-variable relationships. We also: (1) discuss the association between the K-plot and the ROC curve methodology, an aspect that enables us to propose a novel index for measuring dependence; (2) derive the mathematical properties of the proposed measure of dependence.

In this paper, our emphasis has been on constructing a nonparametric measure of dependence based on Kendall's approach. In general, in this framework, there are not most powerful decision making mechanisms. For example, tests such as the Spearman rank-correlation, Kendall-tau, and Fisher-Yates normal scores tests can provide relatively low power levels, when, e.g., $Y$ has a non-monotone regression on $X$ (e.g., Feuerverger, 1993; Vexler et al., 2014). The use of the ROC/AUC methodology in our development may require careful considerations of the proposed method when underlying data correspond to potentially non-monotonic complex dependence structures. However, the proposed technique can provide efficient outputs detecting non-monotonic relationships between the components of a bivariate random variable (for example, see Table 2h).

Note that, in situations when the classical measures of dependence or MIC provide powerful outputs, the proposed approach is still meaningful as a graphical concept incorporating information regarding dependence structures into data analysis.

To demonstrate the applicability of the proposed method we apply our method to the study of biomarkers for myocardial infarction.

The **D**-vector based approach can also be adopted as a valuable tool in the development of linear combinations of biomarkers to maximize their association with a disease factor that follows a continuous distribution. Linear combinations of biomarkers based on the $\boldsymbol{D}$-vector can be constructed in a similar manner to that used in the ROC curve methodology for maximizing the AUC (Chen et al., 2015). In this context, one can consider AUKs based on combinations of random variables to represent associations between groups of variables. Furthermore, the proposed approach can be easily extended to measure dependence of $m$ random variables, say $Z_1, \ldots Z_m$, from a joint distribution function $H$. In this case, AUK $= \Pr\{H(Z_1, \ldots Z_m) < U_1 U_2 \cdots U_m\}$-type objects, where $U_i$, $i = 1, \ldots, m$ are independent uniformly $[0, 1]$ distributed random variables, can be involved. Partial AUCs are known elements in the ROC curve analysis. It may be of practical interest to define and examine partial AUKs, for example, focusing on segments of random variables' distributions, e.g. tails. Further studies are needed to evaluate the proposed approach in these frameworks.

**Appendix A  Proofs**

*Proof of Example included in Section 3.* According to the example's statement, we consider the bivariate random variable $(X, Y = 1 - |X|)$, where $X \sim Unif[-1, 1]$. In this case, we have

$$
\begin{aligned}
H(x, y) &= \Pr\{X < x, Y < y\} = \Pr\{X < x, X > 1 - y, X > 0\} + \Pr\{X < x, X < y - 1, X \leq 0\} \\
&= \Pr\{1 - y < X < x\} I\{x > 0\} + \Pr\{X < \min(x, y - 1)\} \\
&= \frac{1}{2}(x + y - 1) I\{x > 1 - y\} + \frac{1}{2}(\min(x, y - 1) + 1),
\end{aligned}
$$

for $-1 \leq x \leq 1$ and $0 \leq y \leq 1$, where $I\{.\}$ is the indicator function. Then

$$
\begin{aligned}
K(t) &= \Pr\{H(X, Y) < t\} = \Pr\left\{\frac{1}{2}(-|X| + 1) < t\right\} = 1 - \Pr\{2t - 1 < X < 1 - 2t\} \\
&= I\{t > 1/2\} + 2t I\{t \leq 1/2\}, \ 0 \leq t \leq 1.
\end{aligned}
$$

$\square$

*Proof of Proposition 1.* By virtue of the definition of the AUK, we have

$$
\begin{aligned}
\text{AUK} &= \int_0^1 \Pr\{H(X, Y) < t\} \mathrm{d} \Pr(VU < t) = 1 - \Pr\{VU < H(X, Y)\} \\
&= 1 - \mathrm{E}[\Pr\{VU < H(X, Y) | H(X, Y)\}] = 1 - \mathrm{E}[H(X, Y) - H(X, Y) \log\{H(X, Y)\}], \square
\end{aligned}
$$

since $\Pr\{VU < t\} = t - t \log(t)$, for a fixed $t$. The proof is completed.

*Proof of Proposition 2.* By virtue of the Fréchet-Hoeffding upper-bound, we have

$$\Pr\{H(X,Y) < t\} \geq \Pr[\min\{F(X), G(Y)\} < t] \geq \Pr\{F(X) < t\} = t.$$

It is clear that $\Pr\{H(X,Y) < t\} \leq 1$. Then,

$$0 \leq \int_0^1 \Pr\{H(X,Y) \geq t\}\mathrm{d}(t - t\log(t)) \leq \int_0^1 (1-t)(-\log(t))\mathrm{d}t = 3/4.$$

Thus, the range of AUK is $[1/4, 1]$.

In the case $(X, Y = X)$, we have $H(x,y) = \Pr(X < x, X < y)$ and $H(X,Y) = H(X,X)$. Then, $\Pr\{H(X,Y) > t\} = \Pr\{F(X) > t\}$ and $\int_0^1 \Pr\{H(X,Y) \geq t\}\mathrm{d}(t - t\log(t)) = \int_0^1 (1-t)(-\log(t))\mathrm{d}t = 3/4$. Hence, AUK $= 1/4$.

In the case $(X, Y = -X)$, we have $H(x,y) = \Pr(X < x, -X < y)$ and $H(X,Y) = H(X,-X)$, which leads to $\Pr\{H(X,Y) > t\} = 0$ and AUK $= 1$.

Applying these results to the definitions of $\mathrm{AUK}_i$, $i = 1, 2, 3$, we complete the proof. □

*Proof of Proposition 3.* We have

$$\begin{aligned}
H_{12,0}(z_1 + z_2, w_1 + w_2) &= \Pr(Z_1 + Z_2 < z_1 + z_2, W_1 + W_2 < w_1 + w_2) \\
&\geq \Pr(Z_1 + Z_2 < z_1 + z_2, W_1 + W_2 < w_1 + w_2, Z_1 < z_1, W_1 < w_1) \\
&\geq \Pr(Z_1 + Z_2 < Z_1 + z_2, W_1 + W_2 < W_1 + w_2, Z_1 < z_1, W_1 < w_1) \\
&= \Pr(Z_2 < z_2, W_2 < w_2, Z_1 < z_1, W_1 < w_1) \\
&= H_{1,0}(z_1, w_1) H_{2,0}(z_2, w_2).
\end{aligned}$$

Write $H_{12,0} = H_{12,0}(Z_1 + Z_2, W_1 + W_2)$ and $H_{i,0} = H_{i,0}(Z_i, W_i)$, $i = 1, 2$. The preceding inequality and Proposition 1 imply

$$\begin{aligned}
\mathrm{AUK}_{12,0} &= \mathrm{E}\{1 - H_{12,0} + H_{12,0}\log(H_{12,0})\} \\
&\leq \mathrm{E}[1 - H_{1,0}H_{2,0} + H_{1,0}H_{2,0}\{\log(H_{1,0}) + \log(H_{2,0})\}],
\end{aligned}$$

since the function $f(H) = 1 - H + H\log(H)$ is monotone decreasing for $H \in [0,1]$. Then

$$\begin{aligned}
\text{AUK}_{12,0} &\leq \text{E}\{1 - H_{1,0}H_{2,0} + H_{1,0}\log(H_{1,0}) + H_{2,0}\log(H_{2,0})\} \\
&= \text{E}\{1 - H_{1,0} + (H_{1,0} - H_{1,0}H_{2,0}) + H_{1,0}\log(H_{1,0}) + H_{2,0}\log(H_{2,0})\} \\
&= \text{E}\{1 - H_{1,0} + H_{1,0}(1 - H_{2,0}) + H_{1,0}\log(H_{1,0}) + H_{2,0}\log(H_{2,0})\} \\
&\leq \text{E}\{1 - H_{1,0} + (1 - H_{2,0}) + H_{1,0}\log(H_{1,0}) + H_{2,0}\log(H_{2,0})\} \\
&= \text{AUK}_{1,0} + \text{AUK}_{2,0}.
\end{aligned}$$

The proof is complete. □

*Proof of Proposition 4.* Taking into account the statement of requirement $R_6$, we define the transformation $(X^*, Y^*) = (f(X), g(Y))$ and the joint distribution function

$$H^*(x^*, y^*) = \Pr(X^* \leq x^*, Y^* \leq y^*),$$

where $f$ and $g$ are strictly monotone functions on Range $X$ and Range $Y$, respectively. In a similar manner to the $K_i(t), i = 0, 1, 2, 3$, definitions shown in Section 3, we denote the Kendall distribution type functions based on $(X^*, Y^*)$ in the form $K_i^*(t), i = 0, 1, 2, 3$. Consider the following scenarios.

(i) $f$ and $g$ are strictly increasing functions. Then

$$H^*(x^*, y^*) = \Pr\{f(X) \leq f(x), g(Y) \leq g(y)\} = \Pr(X \leq x, Y \leq y) = H(x, y).$$

Thus, $K_0^*(t) = K_0(t)$, $t \in [0,1]$. Similarly, one can easily show that $K_i^*(t) = K_i(t)$, $t \in [0,1]$, for all $i = 1, 2, 3$.

(ii) $f$ is a strictly increasing function, $g$ is a strictly decreasing function. In a similar manner to the case considered above, we conclude that $K_0^*(t) = K_2(t)$, $K_1^*(t) = K_3(t)$, $K_2^*(t) = K_0(t)$ and $K_3^*(t) = K_1(t)$, $t \in [0,1]$.

(iii) $f$ is a strictly decreasing function, $g$ is a strictly increasing function. It is clear that, in this case, $K_0^*(t) = K_1(t)$, $K_1^*(t) = K_0(t)$, $K_2^*(t) = K_3(t)$ and $K_3^*(t) = K_2(t)$, $t \in [0,1]$.

(iv) $f$ and $g$ are strictly decreasing functions. Then $K_0^*(t) = K_3(t)$, $K_1^*(t) = K_2(t)$, $K_2^*(t) = K_1(t)$ and $K_3^*(t) = K_0(t)$, $t \in [0,1]$.

Thus, by virtue of the definition of $I_{\text{AUK}}$ based on the vector $\boldsymbol{D}$, we have that $I_{\text{AUK}}$ satisfies $R_6$.

Regarding $R_4$, we assume the condition $C_1$ is in effect. The $\text{AUK}_i, i = 0, 1, 2, 3$, were designed to be 1/2 when $X$ and $Y$ are independent. Then, in this case, $I_{\text{AUK}} = 0$. Suppose we observe $I_{\text{AUK}} = 0$. Note that $I_{\text{AUK}} = 0$ if and only if $\text{AUK}_i = 1/2$, for all $i = 0, 1, 2, 3$, since $I_{\text{AUK}}^2 = (8/5) \sum_{i=0}^{3} (\text{AUK}_i - 1/2)^2$. Thus, $\int_0^1 K_{i_0}(t) \, \mathrm{d}\,(t - t\log(t)) = 1/2$, whereas it is assumed that $K_{i_0}(t) \geq (\text{or } \leq) \{t - t\log(t)\}$, for all $0 \leq t \leq 1$. Then $K_{i_0}(t) = \{t - t\log(t)\}$, for all $0 \leq t \leq 1$. Therefore, the components of $K_{i_0}$ related random pair $((X,Y)$ or $(X,-Y)$, or $(-X,Y)$, or $(-X,-Y))$ are independent, which implies that $X$ and $Y$ are independent.

Regarding R$_8$, we note that by virtue of the weak convergence of $H_n$ to $H$, we have $H_{i,n}$ converges weakly to $H_i$, $i = 0, \ldots, 3$, where $H_{i,n}, i = 1 \ldots, 3$, correspond to the functions $H_i, i = 1 \ldots, 3$, defined in Section 3. Since, for any $\delta > 0$,

$$|1 - H_{i,n}(X_n, Y_n) - H_{i,n}(X_n, Y_n) \log\{H_{i,n}(X_n, Y_n)\}|^{1+\delta} \leq 1,$$

the sequence $1 - H_{i,n}(X_n, Y_n) - H_{i,n}(X_n, Y_n) \log\{H_{i,n}(X_n, Y_n)\}$ is uniformly integrable. Then, using Proposition 1, we define

$$\text{AUK}_{i,n} = \mathrm{E}[1 - H_{i,n}(X_n, Y_n) - H_{i,n}(X_n, Y_n) \log\{H_{i,n}(X_n, Y_n)\}],$$

obtaining that $\lim_{n \to \infty} \text{AUK}_{i,n} = \mathrm{E}[1 - H_i(X,Y) - H_i(X,Y)\log\{H_i(X,Y)\}] = \text{AUK}_i$. The index $I_{\text{AUK}}$ is a continuous function of the corresponding AUK's. Thus, $\lim_{n \to \infty} I_{\text{AUK}}(X_n, Y_n) = I_{\text{AUK}}$. Requirement R$_8$ is satisfied. This completes the proof. $\square$

**Acknowledgement**

Dr. Vexler's effort was supported by the National Institutes of Health (NIH) grant 1G13LM012241-01. Dr. Markatou would like to thank the Jacobs School of Medicine and Biomedical Sciences for facilitating this work through institutional financial resources (to M. Markatou) that supported the work of the second author of this paper. She also acknowledges PCORI for partially supporting this work under award 1507-31640. The authors are grateful to the Editor, the Associate Editor and the referees for suggestions that led to a substantial improvement of this paper.

**Supporting information**

Additional supporting material may be found in the online version of this article:

***Mathematica* and R *Codes:*** Software codes to implement the developed method in this article.

# Supplementary Material:
# Multi-Panel Kendall Plot in Light of an ROC Curve Analysis Applied to Measuring Dependence

Albert Vexler*

Department of Biostatistics,The State University of New York at Buffalo, USA

*Corresponding author. Email: avexler@buffalo.edu

Georgios Afendras

Department of Biostatistics and Jacobs School of
Medicine and Biomedical Sciences,The State University of New York at Buffalo, USA

and

Marianthi Markatou

Department of Biostatistics and Jacobs School of
Medicine and Biomedical Sciences,The State University of New York at Buffalo, USA

**Abstract**

This online supplement to "Multi-Panel Kendall Plot in Light of an ROC Curve Analysis Applied to Measuring Dependence" contains Mathematica (Wolfram, 1999) and *R* (R Development Core Team, 2008) codes to implement the developed method that is proposed in the article.

## SM.1  MATHEMATICA AND *R* CODES

### SM.1.1  The Mathematica Code used in Section 4

In[1]:= $F[u1\_, u2\_, \gamma\_] = (1 - (u1 * u2 * (1 + \gamma * (1 - u1) * (1 - u2)))) + (u1 * u2 * (1 + \gamma * (1 - u1) * (1 - u2)))$
$\text{Log}[(u1 * u2 * (1 + \gamma * (1 - u1) * (1 - u2)))]) * (1 + \gamma * (1 - 2 * u1) * (1 - 2 * u2))$
$\text{Integrate}[F[u1, u2, \gamma], \{u1, 0, 1\}, \{u2, 0, 1\}]$

Out[1]:= $\frac{1}{216\gamma^2}\left(\gamma(-156 + (501 - 58\gamma)\gamma) - 6\pi^2(1 + 3(-2 + \gamma)\gamma) + 12((1 + \gamma)(10 + (-19 + \gamma)\gamma)\right.$
$\left. + 3(1 + 3(-2 + \gamma)\gamma)\text{Log}[-\gamma]\text{Log}[1 + \gamma] + 36(1 + 3(-2 + \gamma)\gamma)\text{PolyLog}[2, 1 + \gamma]\right)$

In[2]:= $\text{Apart}[\text{Integrate}[F[u1, u2, \gamma], \{u1, 0, 1\}, \{u2, 0, 1\}], \gamma]$

Out[2]:= $\frac{1}{108}\gamma(-29 + 6\text{Log}[1 + \gamma]) + \frac{-13 + 3\pi^2 - 9\text{Log}[1 + \gamma] - 18\text{Log}[-\gamma]\text{Log}[1 + \gamma] - 18\text{PolyLog}[2, 1 + \gamma]}{18\gamma}$
$+ \frac{-\pi^2 + 20\text{Log}[1 + \gamma] + 6\text{Log}[-\gamma]\text{Log}[1 + \gamma] + 6\text{PolyLog}[2, 1 + \gamma]}{36\gamma^2}$
$+ \frac{1}{72}(167 - 6\pi^2 - 72\text{Log}[1 + \gamma] + 36\text{Log}[-\gamma]\text{Log}[1 + \gamma] + 36\text{PolyLog}[2, 1 + \gamma])$

In[3]:= $\text{Plot}[\%2, \{\gamma, -0.99, 0.99\}]$

Out[3]:=

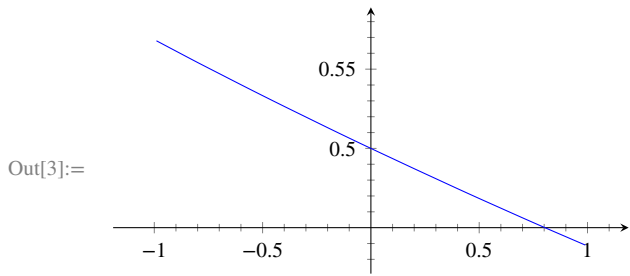

## SM.1.2 The *R* Code used in Section 4

```
AUK<-function(a)
    {
    C<-function(u1,u2) u1*u2*(1+a*(1-u1)*(1-u2))
    G0<-function(u1,u2) (1-C(u1,u2)+C(u1,u2)*log(C(u1,u2)))*(1+a*(1-2*u1)*(1-2*u2))
    G1<-function(u1)
        {
        G2<-function(u2) Vectorize(G0(u1,u2))
        R<-integrate(G2,0,1)$value
        return(R)
        }
    G1<-Vectorize(G1)
    RR<-integrate(G1,0,1)$value
    return(RR)
    }
AUK<-Vectorize(AUK)
A<-seq(from=-1,to=1,by=0.005)
AUKV<-array()
AUKV<-AUK(A)
plot(A,AUKV)
```

### The Monte Carlo Method

```
N<-50000
U1<-runif(N)
U2<-runif(N)
AUK<-function(a)
    {
    C<-function(u1,u2) u1*u2*(1+a*(1-u1)*(1-u2))
    G0<-function(u1,u2) (1-C(u1,u2)+C(u1,u2)*log(C(u1,u2)))*(1+a*(1-2*u1)*(1-2*u2))
    G1<-Vectorize(G0)
    RR<-mean(G1(U1,U2))
    return(RR)
    }
AUK<-Vectorize(AUK)
A<-seq(from=-1,to=1,by=0.005)
AUKV<-array()
AUKV<-AUK(A)
plot(A,AUKV)
AUKV
```

## SM.1.3 The *R* Code Example to Compute AUK Based on $N_2\left(\begin{pmatrix}0\\0\end{pmatrix},\begin{pmatrix}1 & \rho\\ \rho & 1\end{pmatrix}\right)$ Distribution

```
library(pbivnorm)
library(mvtnorm)
g<- function(x){ifelse(x<=0,0,1-x+x*log(x))}
g<- Vectorize(g)
#
AUK<- function(rho){
  Sigma<- matrix(c(1,rho,rho,1), 2, 2)
  phi<- function(x,y){dmvnorm(c(x,y), mean = c(0, 0), sigma = Sigma, log = FALSE)}
  phi<- Vectorize(phi)
  H<- function(x,y){pbivnorm(x,y, rho)}
  H<- Vectorize(H)
auk<- integrate(function(y){
              sapply(y, function(y){integrate(function(x)
              g(H(x,y))*phi(x,y),-Inf,Inf)$value})},-Inf,Inf)$value
return(auk)
}
```

```
#
AUK<- Vectorize(AUK)
```

## The Monte Carlo Method

```
library(MASS)
library(mvtnorm)
#
p<-(-0.5)
Sigma <- matrix(c(1,p,p,1),2,2)
#
N1<-30000
N2<-5000
#
X<-mvrnorm(N1, rep(0, 2), Sigma)
x<-X[,1]
y<-X[,2]
#
C1<-function(u1,u2)  mean(1*((x<u1)&(y<u2))) #nonparametric estimation of H
C11<-Vectorize(C1)
#
X0<-mvrnorm(N2, rep(0, 2), Sigma)
x0<-X0[,1]
y0<-X0[,2]
#
G<-array()
G<-C11(x0,y0)
#
GlogG<-array(0,N2)
GlogG[G!=0]<-G[G!=0]*log(G[G!=0])
MF<-mean(1-G+GlogG)
print(MF)
```